\documentclass[11pt]{article}

\usepackage[small,bf]{caption}
\usepackage{amsmath, amsthm, amssymb}
\usepackage{algorithmic}
\usepackage{algorithm}
\usepackage[parfill]{parskip}

\usepackage{epsfig,fullpage,psfrag,url,verbatim}
\usepackage{graphics}
\usepackage{graphicx}

\usepackage{amsfonts}
\usepackage{color}

\newtheorem{thm}{Theorem}[section]
\newtheorem{lem}{Lemma}[section]
\newtheorem{prop}{Proposition}[section]

\newtheorem{mydef}{Definition}[section]
\newtheorem{rem}{Remark}[section]

\newcommand{\diag}{\mathrm{diag}}
\newcommand{\BDiag}{\mathrm{{\bf D}iag}}
\newcommand{\BMdiag}{\mathrm{{\bf M}diag}}

\newcommand{\hh}{h(\cdot)}
\newcommand{\ff}{f(\cdot)}
\newcommand{\EE}{\mathbb{E}}
\newcommand{\hx}{\hat{x}}
\newcommand{\Dh}{D_h}
\newcommand{\intQ}{\mathrm{int}\ Q}

\newcommand{\bbf}{\bar{f}}
\newcommand{\bh}{\bar{h}}

\newcommand{\bbR}{\mathbb{R}}

\DeclareMathOperator*{\argmin}{arg\,min}

\begin{document}

\title{Relatively-Smooth Convex Optimization by First-Order Methods, and Applications}

\author{Haihao Lu\thanks{MIT Department of Mathematics, 77 Massachusetts Avenue, Cambridge, MA   02139
({mailto:  haihao@mit.edu}).}
\and
Robert M. Freund\thanks{MIT Sloan School of Management, 77 Massachusetts Avenue, Cambridge, MA   02139
({mailto:  rfreund@mit.edu}).  This author's research is supported by AFOSR Grant No. FA9550-15-1-0276 and  the MIT-Belgium Universit\'{e} Catholique de Louvain Fund.}
\and
Yurii Nesterov\thanks{Universit\'{e} Catholique de Louvain
({mailto:  yurii.nesterov@uclouvain.be}).  This author's research is supported by the MIT-Belgium Universit\'{e} Catholique de Louvain Fund.}}
\date{Revised October 7, 2017 (original dated October 17, 2016)} 

\maketitle

\begin{abstract}

The usual approach to developing and analyzing first-order methods for smooth convex optimization assumes that the gradient of the objective function is uniformly smooth with some Lipschitz constant $L$. However, in many settings the differentiable convex function $f(\cdot)$ is not uniformly smooth -- for example in $D$-optimal design where $f(x):=-\ln \det(HXH^T)$ and $X:=\BDiag (x)$, or even the univariate setting with $f(x) := -\ln(x) + x^2$. In this paper we develop a notion of ``relative smoothness'' and relative strong convexity that is determined relative to a user-specified ``reference function'' $h(\cdot)$ (that should be computationally tractable for algorithms), and we show that many differentiable convex functions are relatively smooth with respect to a correspondingly fairly-simple reference function $h(\cdot)$.  We extend two standard algorithms -- the primal gradient scheme and the dual averaging scheme -- to our new setting, with associated computational guarantees.  We apply our new approach to develop a new first-order method for the $D$-optimal design problem, with associated computational complexity analysis.  Some of our results have a certain overlap with the recent work \cite{bbt}.

\end{abstract}

\section{Introduction, Definition of ``Relative-Smoothness,'' and Basic Properties}\label{intro}

\subsection{Traditional Set-up for Smooth First-Order Methods}\label{traditional}
Our optimization problem of interest is:
\begin{equation}\label{poi1}
\begin{array}{lrlr} P:  & f^* := \ \  \mbox{minimum}_x & f(x) \\ \\
& \mbox{ s.t. } & x \in Q \ ,
\end{array}
\end{equation}
where $Q \subseteq \mathbb{E}$ is a closed convex set in the finite-dimensional vector space $\mathbb{E}$  with inner product $\langle \cdot\  ,\  \cdot \rangle$, and $f(\cdot) : Q \to \mathbb{R}$ is a differentiable convex function.

There are by now very many first-order methods for tackling the optimization problem \eqref{poi1}, see for example \cite{nesterovBook}, \cite{tseng}, \cite{polyak1987introduction};  virtually all such methods are designed to solve \eqref{poi1} when the gradient of $f(\cdot)$ satisfies a uniform Lipschitz condition on $Q$, namely there exists a constant $L_f < \infty$ for which:
\begin{equation}\label{lipschitz}
\|\nabla f(x) - \nabla f(y)\|_* \le L_f \|x-y\| \ \ \ \text{for~all~} x,y \in Q \ ,
\end{equation} where $\|\cdot\|$ is a given norm on $\mathbb{E}$ and $\|\cdot\|_*$ denotes the usual dual norm.  For example, consider the standard gradient descent scheme, which presumes the norm in \eqref{lipschitz} is Euclidean, and uses the following update:
\begin{equation}
x^{i+1} \gets \arg\min_{x \in Q} \left\{f(x^i) + \langle \nabla f(x^i), x-x^i \rangle+ \tfrac{L_f}{2}\|x-x^i\|_2^2 \right\} \ .
\end{equation} One can prove for the standard gradient descent scheme that after $k$ iterations it holds for any $x \in Q$ that:
\begin{equation}\label{sd_bound3}
 f(x^k)  - f(x) \ \ \leq \ \  \frac{L_f \|x-x^0\|_2^2}{2k} \ ,
\end{equation} which is an $O(1/k)$ sublinear rate of convergence \cite{nesterovBook}, \cite{polyak1987introduction}.  Furthermore, if $f(\cdot)$ is also uniformly $\mu_f$-strongly convex for some $\mu_f >0$, namely:
\begin{equation}\label{arm}
f(y) \ge f(x) + \langle \nabla f(x), y-x \rangle + \tfrac{\mu_f}{2} \|y-x\|^2_2 \   \ \ \ \text{for~all~} x,y \in Q \ ,
\end{equation}then one can prove linear convergence for the gradient descent scheme, see \cite{nesterovBook}, \cite{polyak1987introduction}, i.e., for any $x \in Q$ we have that:
\begin{equation}\label{sd_bound4}
f(x^k)  - f(x) \ \ \leq \ \  \frac{L_f}{2}\left(1-\frac{2\mu_f}{L_f+\mu_f}\right)^{k} \|x-x^0\|_2^2 \ .
\end{equation}

More general versions of first-order methods are not restricted to the Euclidean ($\|\cdot\|_2$) norm, and use a differentiable ``prox function'' $h(\cdot)$, which is a $1$-strongly convex function on $Q$, to define a Bregman distance:
\begin{equation}\label{bregman}
\Dh(y,x) := h(y) - h(x) - \langle \nabla h(x), y-x \rangle \ \ \ \text{for~all~} x,y \in Q
\end{equation} which as a result satisfies $$\Dh(y,x) \ge \tfrac{1}{2} \|y-x\|^2 \ . $$ The standard Primal Gradient Scheme (with Bregman distance), see \cite{tseng}, has the following update formula:
\begin{equation}\label{proxgrad}
x^{i+1} \gets \arg\min_{x \in Q} \left\{f(x^i) + \langle \nabla f(x^i) , x-x^i \rangle + L_f \Dh(x,x^i) \right\} \ .
\end{equation}
Notice in \eqref{proxgrad} by construction that the update requires the capability to solve instances of a subproblem of the general form:
\begin{equation}\label{subproblem}
x_{\text{new}} \gets \arg\min_{x \in Q}\{\langle c, x \rangle + h(x)\} \ ,
\end{equation} for suitable iteration-specific values of $c$; indeed, \eqref{proxgrad} is an instance of the subproblem \eqref{subproblem} with $c = \tfrac{1}{L_f}\nabla f(x^i) - \nabla h(x^i)$ at iteration $i$.  It is especially important to note that the Primal Gradient Scheme is somewhat meaningless whenever we do not have the capability to efficiently solve \eqref{subproblem}, a point which we will return to later on.  In a typical design and implementation of a first-order method for solving \eqref{poi1}, one attempts to specify the norm $\|\cdot\|$ and the strongly convex prox function $h(\cdot)$ in consideration of the shape of the feasible domain $Q$ while also ensuring that the subproblem \eqref{subproblem} is efficiently solvable.

Regarding computational guarantees, one can prove for the Primal Gradient Scheme that after $k$ iterations it holds for any $x \in Q$ that:
\begin{equation}\label{sd_bound33}
 f(x^k)  - f(x) \ \ \leq \ \  \frac{L_f \Dh(x,x^0)}{k} \ ,
\end{equation} which is an exact generalization of \eqref{sd_bound3}, see \cite{tseng}, \cite{NemirovskyYudin83}.

We emphasize that standard first-order methods as stated above for solving \eqref{poi1} require that $f(\cdot)$ be uniformly smooth on $Q$, that is, that there is a finite value of the Lipschitz constant $L_f$ as defined in \eqref{lipschitz}, in order to ensure associated computational guarantees. However, there are many differentiable convex functions in practice that do not satisfy a uniform smoothness condition.  Consider $f(x):=-\ln \det(HXH^T)$ with $X:=\BDiag(x)$ in $D$-optimal design on the feasible set $Q = \{ x \in \bbR^n : \langle e, x \rangle = 1, \ x \ge 0 \}$, or $f(x)= |x|^3$ or $f(x) = x^4$ on the feasible set $Q=\bbR$, or $f(x) = -\ln(x) + x^2$ on $Q=\bbR_{++}$.  Of course, if the algorithm iterates have monotone decreasing objective function values (which is provably the case for most smooth first-order methods), it then is sufficient just to ensure that $f(\cdot)$ is smooth on some level set of $f(\cdot)$.  Nevertheless, even in this case the constant $L_f$ may be huge.  For instance, let  $f(x) = -\ln(x) + x^2$ on $Q=\bbR_{++}$, and consider the level set $\{x : f(x)\le 10\}$. Then one still has $L_f\approx \exp^{20}$ on this level set, which is not reasonable for practical use.

Notice that unlike quadratic functions, the second-order terms of the functions in the above examples vary dramatically on $Q$ -- and especially as $x \rightarrow \partial Q$ (or as $x$ goes to infinity in $Q$). It therefore becomes unreasonable to use a uniform bound of the form $L_f$ to upper-bound second-order information.  

Motivated by the above drawbacks in standard first-order methods, we develop a notion of ``relative smoothness'' and relative strong convexity, relative to a given ``reference function'' $h(\cdot)$ and which does not require the specification of any particular norm -- and indeed $h(\cdot)$ need not be either strictly or strongly convex.  Armed with relative smoothness and relative strong convexity, we demonstrate the capability to solve a more general class of differentiable convex optimization problems (without uniform Lipschitz continuous gradients), and we also demonstrate linear convergence results for both a Primal Gradient Scheme and a Dual Averaging Scheme when the function is both relatively smooth and relatively strongly convex.

There is a certain overlap of ideas and results herein with the paper \cite{bbt} by Bolte, Bauschke, and Teboulle.  For starters, the relative smoothness condition definition in the present paper in Definition \ref{smooth} is equivalent to the (LC) condition in \cite{bbt} except that \cite{bbt} also requires the reference function $\hh$ to be essentially smooth and strictly convex, which we do not need in this paper.  The main developments in \cite{bbt} are based on generalizing a key descent lemma and applying this generalization to tackle (additive) composite optimization problems using the primal gradient scheme (called the NoLips Algorithm in \cite{bbt}) with associated complexity analysis involving a symmetry measure of the Bregman distance $D_h(\cdot, \cdot)$.  These results are then illustrated in the application of composite optimization to Poisson inverse problems.  While the NoLips Algorithm in \cite{bbt} is structurally the same as Algorithm \ref{proxgradscheme} herein, they are both instantiations of the standard primal gradient scheme; however, as will be seen in Section \ref{algorithms} here, we do not need any symmetry measure in constructing step-sizes or in the complexity analysis.  The paper \cite{zhouaunified} by Zhou, Liang, and Shen also tackles composite optimization using the standard primal gradient scheme which therein is called PGA-$\cal B$, with a focus on demonstrating equivalence of proximal gradient and proximal point methods more broadly. Here we develop  measures of relative smoothness and also relative strong convexity, which can improve the computational guarantees of the primal gradient scheme, see Theorem \ref{thmprimalgd}.  We further present computational guarantees for the dual averaging scheme \cite{nesterovda} in Theorem \ref{thmdualgd}.  In Section \ref{examples} we show that many differentiable convex functions are relatively smooth with respect to a correspondingly fairly-simple reference function $h(\cdot)$ that is easy to construct and for which algorithmic computations can be effeciently be performed.  In Section \ref{return} we apply our approach to develop a new first-order method for the $D$-optimal design problem, with associated computational complexity analysis.  Throughout the current paper, we compare and clarify similarities and differences between our work and \cite{bbt} in the context of the specific contributions as they arise.

\subsection{Relative Smoothness and Relative Strong Convexity}



Let $\hh$ be any given differentiable convex function (it need not be strongly nor even strictly convex) defined on $Q$.  We will henceforth refer to $\hh$ as the ``reference function.''  We define ``relative smoothness'' and ``relative strong convexity'' of $f(\cdot)$ relative to $\hh$ using the Bregman distance \eqref{bregman} associated with $\hh$ as follows.\medskip\medskip

\begin{mydef}\label{smooth} $\ff$ is $L$-smooth relative to $\hh$ on $Q$ if for any
$x,y\in \intQ$, there is a scalar $L$ for which
\begin{equation}\label{robocop}
f(y) \le f(x) + \langle \nabla f(x), y-x \rangle +  L \Dh(y,x) \ .
\end{equation}
\end{mydef}\medskip

\begin{mydef}\label{strong}
$\ff$ is $\mu$-strongly convex relative to $\hh$ on $Q$ if for any
$x,y\in \intQ$, there is a scalar $\mu  \ge 0$ for which
\begin{equation}
f(y) \ge f(x) + \langle \nabla f(x), y-x \rangle +  \mu \Dh(y,x) \ .
\end{equation}
\end{mydef}

(Here and elsewhere $\intQ$ denotes the interior of $Q$.  In cases where $Q$ has no interior, one can instead use the relative interior of $Q$.)  Note that relative smoothness and relative strong convexity of $\ff$ are defined relative to the reference function $\hh$ directly; no norm is involved in the definitions, so that smoothness/strong convexity does not depend on any norm.  Furthermore, $\hh$ is not presumed to have any special properties by itself such as strict or (traditional) strong convexity; rather the key structural properties involve how $\ff$ behaves relative to $\hh$.  The definition of relative smoothness above is equivalent to the (LC) condition in \cite{bbt}, but \cite{bbt} requires the reference function to be essentially smooth and strictly convex, which we do not need.

The following proposition presents equivalent definitions of relative smoothness and relative strong convexity.  In the case when both $\ff$ and $\hh$ are twice differentiable, parts (a-iii) and (b-iii) of the proposition demonstrate that the above definitions are equivalent to $$\mu \nabla^2 h(x) \preceq  \nabla^2 f(x) \preceq L \nabla^2 h(x) \ \ \text{for~all} \ x \in \intQ \ , $$ which is an intuitively simple condition on the Hessian matrices of the two functions.\medskip

\begin{prop}\label{thm:equivdef}
\begin{description}
\item[ ] The following conditions are equivalent:
\begin{description}
\item[(a-i)] $\ff$ is $L$-smooth relative to $\hh$,
\item[(a-ii)] $L\hh - \ff$ is a convex function on $Q$,
\item[(a-iii)] Under twice-differentiability $\nabla^2 f(x) \preceq L\nabla^2 h(x)$ for any $x\in \intQ$,
\item[(a-iv)] $\langle \nabla f(x)-\nabla f(y), x-y \rangle \le L \langle \nabla h(x) - \nabla h(y), x-y \rangle$ for all $x,y \in \intQ$.
\end{description}
\item[ ] The following conditions are equivalent:
\begin{description}
\item[(b-i)] $\ff$ is $\mu$-strongly convex relative to $\hh$,
\item[(b-ii)] $\ff - \mu \hh$ is a convex function on $Q$,
\item[(b-iii)]  Under twice-differentiability $\nabla^2 f(x) \succeq \mu\nabla^2 h(x)$ for any $x\in \intQ$,
\item[(b-iv)] $\langle \nabla f(x)-\nabla f(y), x-y \rangle \ge \mu \langle \nabla h(x) - \nabla h(y), x-y \rangle$ for all $x,y \in \intQ$. \qed
\end{description}
\end{description}
\end{prop}\medskip\medskip

The first part of Proposition \ref{thm:equivdef} is almost equivalent to Proposition 1 of \cite{bbt}.

{\bf Proof:}  For $x \in Q$ define $\phi(x):= Lh(x) - f(x)$.  Using \eqref{robocop} and \eqref{bregman} it follows that (a-i) holds if and only if $\phi(x) \ge \phi(y) + \langle \nabla \phi(y),x-y\rangle $ for all $x,y \in Q$, which is equivalent to the convexity of $\phi(\cdot) = Lh(\cdot) - f(\cdot)$ from Theorem 2.1.2 of \cite{nesterovBook}, thus showing that (a-i) $\Leftrightarrow$ (a-ii).  It follows from Theorem 2.1.3 of \cite{nesterovBook} applied to $\phi(\cdot)$ that $\phi(\cdot)$ is convex if and only if   $ \langle \nabla \phi(x) - \nabla \phi(y),x-y\rangle \ge 0 $ for all $x,y \in Q$, which shows that (a-ii) $\Leftrightarrow$ (a-iv).  If $f(\cdot)$ and $h(\cdot)$ are twice differentiable, then it follows from Theorem 2.1.4 of \cite{nesterovBook} that (a-ii) $\Leftrightarrow$ (a-iii).

Similar proofs can be applied for part (b).\qed

%
%

For notational convenience, let us denote by $\ff \preceq \hh$ that $\hh - \ff$ is a convex function, whereby this also means $\ff$ is $1$-smooth with respect to $\hh$ from Proposition \ref{thm:equivdef}.  Similarly $\ff \succeq \hh$ means $\ff - \hh$ is a convex function and so $\ff$ is $1$-strongly convex with respect to $\hh$. (In the case when both $\ff$ and $\hh$ are twice differentiable, the relation ``$ \cdot \succeq \cdot$'' on two functions is consistent with the L\"{o}wner partial order on the Hessians of these two functions from Propositon \ref{thm:equivdef}.)  Then the condition that $\ff$ is $L$-smooth with respect to $\hh$ is equivalent to $\ff \preceq L \hh$; similarly the condition that $\ff$ is $\mu$-strongly convex with respect to $\hh$ is equivalent to $\ff \succeq \mu \hh$. In addition, relative-smoothness and relative strong convexity are each transitive, so that $\ff \preceq g(\cdot)$ and $g(\cdot) \preceq \hh$ implies that $\ff \preceq \hh$.


We can also work with sums and linear transformations of relatively smooth and/or relatively strongly convex functions, as the next proposition states.\medskip

\begin{prop}\label{prop:additivity}
\mbox{}
\begin{description}
\item[1.] If $f_1(\cdot)\preceq L_1\hh$ and $f_2(\cdot) \preceq L_2 h_2(\cdot)$, then for all $\alpha,\beta \ge 0$ it holds that $f(\cdot):=\alpha f_1(\cdot)+\beta f_2(\cdot) \preceq h(\cdot):=\alpha L_1 h_1(\cdot)+\beta L_2 h_2(\cdot)$.
\item[2.] If $f_1(\cdot) \succeq \mu_1h_1(\cdot)$ and $f_2(\cdot) \succeq \mu_2h_2(\cdot)$, then for all $\alpha,\beta \ge 0$ it holds that $f(\cdot):=\alpha f_1(\cdot)+\beta f_2(\cdot) \succeq h(\cdot):=\alpha \mu_1 h_1(\cdot)+\beta \mu_2 h_2(\cdot)$.
\item[3.] If $f(\cdot) \preceq h(\cdot)$, and $A$ is a linear transformation of appropriate dimension, then $\phi_f(x):= f(Ax) \preceq \phi_h(x):=h(Ax)$.
\item[4.] If $f(\cdot)\succeq h(\cdot)$, and $A$ is a linear transformation of appropriate dimension, then $\phi_f(x):= f(Ax) \succeq \phi_h(x):=h(Ax)$.
\end{description}
\end{prop}
{\bf Proof:} The proofs of the first two arguments follow directly from the definitions of relative smoothness and relative strong convexity in Definitions \ref{smooth} and \ref{strong}. The proofs of the last two arguments follow from the equivalent definition (a-iv) and (b-iv) in Proposition \ref{thm:equivdef}.\qed

{

\subsection{Constructive Algorithmic Set-up}

Let us now discuss criteria for choosing the reference function $\hh$ in the context of computational schemes for solving the optimization problem \eqref{poi1}.  To be concrete, consider a simple Primal Gradient Scheme as shown in Algorithm \ref{proxgradscheme}.  Note that this scheme is essentially as described in the update formula \eqref{proxgrad}, except that the uniform smoothness constant $L_f$ is replaced by the relative smoothness parameter $L$ of $\ff$ with respect to the reference function $\hh$ as defined in Definition \ref{smooth}, and the only formal requirement for $\hh$ is that the pair $(\ff, \hh)$ must satisfy the conditions of Definition \ref{smooth}.\medskip

\begin{algorithm}
\caption{Primal Gradient Scheme with reference function $h(\cdot)$}\label{proxgradscheme}
$ \ $
\begin{algorithmic}
\STATE {\bf Initialize.}  Initialize with $x^0 \in Q$.  Let $L$, $\hh$ satisfying Definition \ref{smooth} be given. \\
At iteration $i$ :\\
\STATE  {\bf Perform Updates.}  Compute $\nabla f(x^i)$ , \\\medskip
\ \ \ \ \ \ \ \ \ \ \ \ \ \ \ \ \ \ \ \ \ \ \ \ \ \ \ \ $x^{i+1} \gets \arg\min_{x \in Q} \{f(x^i) + \langle \nabla f(x^i), x-x^i \rangle + L \Dh(x,x^i) \}$ .\\\medskip
\end{algorithmic}
\end{algorithm}\medskip

In order to efficiently execute the update step in Algorithm \ref{proxgradscheme} we also require of $\hh$ that the subproblem \eqref{subproblem} is efficiently solvable for any given $c$.  In summary, to solve the optimization problem \eqref{poi1} using Algorithm \ref{proxgradscheme}, we need to specify a reference function $\hh$ that has the following two properties:
\begin{enumerate}
\item[(i)] $\ff$ is $L$-smooth relative to $\hh$ on $Q$, and
\item[(ii)] the subproblem \eqref{subproblem} always has a solution, and the solution is efficiently computable.
\end{enumerate}

In Section \ref{examples} we will see how this can be done for several useful classes of problems that are not otherwise solvable by traditional first-order methods that require uniform Lipschitz continuity of the gradient.  In Section \ref{algorithms} we analyze the computational guarantees associated with the Primal Gradient Scheme (Algorithm \ref{proxgradscheme}) as well as a Dual Averaging Scheme.  In Section \ref{return}, we apply the computational guarantees of Section \ref{algorithms} to the $D$-optimal design problem.

\noindent {\bf Notation.} For a vector $x$, $X = \BDiag(x)$ denotes the diagonal matrix with the coefficients of $x$ along the diagonal. For a symmetric matrix $A$, $\diag (A)$ denotes the vector of the diagonal coefficients of $A$, and $\BMdiag(A)$ denotes the diagonal matrix whose diagonal coefficients correspond to the diagonal coefficients of $A$.  Unless otherwise specified, the norm of a matrix is the operator norm using $\ell_2$ norms.  The $\ell_p$ norm of a vector $x$ is denoted by $\|x\|_p$.  For symmetric matrices, ``$\succeq$'' denotes the L\"{o}wner partial order. In a mild double use of notation, $\ff \succeq \hh$ denotes $\ff-\hh$ is a convex function, and the appropriate meaning of ``$\succeq$'' will be obvious in context. Let $e$ denote the vector of $1$'s whose dimension is dictated by context. Let $\Delta_n := \{x \in \mathbb{R}^n : \langle e, x \rangle = 1, \ x \ge 0\}$ denote the standard unit simplex in $\mathbb{R}^n$. Given two matrices $A$ and $B$ of the same order, let $ A \circ B$ denote the Hadamard (i.e., component-wise) product of $A$ and $B$, see for example Anstreicher \cite{kurt}.  Let $\exp$ denote the base of the natural logarithm.

\medskip\medskip

\section{Examples of Relatively Smooth Optimization Problems}\label{examples}

Here we show several classes of optimization problems \eqref{poi1} for which one can easily construct a reference function $\hh$ with the two properties mentioned above, namely (i) $\ff$ is $L$-smooth relative to $\hh$ for an easily determined value $L$, and (ii) the subproblem \eqref{subproblem} is efficiently solvable.

\subsection{Optimization over $\mathbb{R}^n$ with $\|\nabla^2f(x)\|$ growing as a polynomial in $\|x\|_2$}\label{sec:rn}

Suppose that $\ff$ is a twice-differentiable convex function on $Q := \mathbb{R}^n$ and let $\| \nabla^2 f(x) \|$ denote the operator norm of $\nabla^2 f(x)$ with respect to the $\ell_2$-norm on $\mathbb{R}^n$.   Suppose that $\| \nabla^2 f(x) \| \le p_r(\|x\|_2)$, where $p_r(\alpha)=\sum_{i=0}^r a_i \alpha^i$ is an $r$-degree polynomial of $\alpha$.  Let \begin{equation}\label{sunday} h(x) :=\tfrac{1}{r+2} \|x\|_2^{r+2} + \tfrac{1}{2}\|x\|_2^2 \ . \end{equation}  Then the following proposition states that $\ff$ is $L$-smooth relative to $\hh$ for an easily computable value $L$.  This implies that no matter how fast the Hessian of $\ff$ grows as $\|x\|_2 \rightarrow \infty$, $\ff$ can still be smooth relative to the simple reference function $\hh$, even though $\nabla f(\cdot)$ need not exhibit uniform Lipschitz continuity. \medskip\medskip

\begin{prop}\label{lem:smoothness-Rn}
Suppose $\ff$ is twice differentiable and satisfies $\| \nabla^2 f(x) \| \le p_r(\|x\|_2)$ where $p_r(\alpha)$ is an $r$-degree polynomial of $\alpha$. Let $L$ be such that $p_r(\alpha)\le L (1+\alpha^r)$ for $\alpha\ge 0$.  Then $\ff$ is $L$-smooth relative to $h(x) =\frac{1}{r+2} \|x\|_2^{r+2} + \frac{1}{2}\|x\|_2^2$.
\end{prop}
{\bf Proof:}
It follows from elementary rules of differentiation that
$$
\nabla^2 h(x) \ = \  (1+\|x\|_2^r) I + (r+1)\|x\|_2^{r-2} x x^T \  \succeq \ (1+\|x\|_2^r) I \  \succeq \ \tfrac{1}{L} p_r (\|x\|_2) I \ \succeq \ \tfrac{1}{L} \nabla^2 f(x) \ ,
$$ and so $\ff$ is $L$-smooth relative to $\hh$ by part (iii) of Proposition \ref{thm:equivdef}.
\qed

Utilizing the additivity property in Proposition \ref{prop:additivity} together with Proposition \ref{lem:smoothness-Rn}, one concludes that virtually every twice-differentiable convex function on $\mathbb{R}^n$ is $L$-smooth relative to some simple polynomial function of $\|x\|_2$.\medskip\medskip

\begin{rem}\label{uno} Suppose $p_r(\alpha)=\sum_{i=0}^r a_i \alpha^i$. In Proposition \ref{lem:smoothness-Rn}, one simple way to set $L$ is to use $L=\sum_{i=0}^{r}|a_i|$.  Then
\begin{equation}
p_r(\alpha)\le \left\{\begin{array}{ll}
\sum_{i=0}^r |a_i| & \mathrm{for} \  0 \le \alpha \le 1 \\ \\
\sum_{i=0}^r |a_i| \alpha^r &  \mathrm{for} \ \alpha \ge 1 \ ,
\end{array}\right.
\end{equation}
whereby $p_r(\alpha) \le \max \{ L, L\alpha^r \} \le L(1+\alpha^r) $ for $\alpha \ge 0$.\end{rem}\medskip

\noindent {\bf Solving the subproblem \eqref{subproblem}.}  Let us see how we can solve the subproblem \eqref{subproblem} for this class of optimization problems.  The subproblem \eqref{subproblem} can be written as
\begin{equation}\label{eq:rnsubprob}
\min_{x \in \mathbb{R}^n} \ \  \langle c, x \rangle + \tfrac{1}{r+2} \|x\|_2^{r+2} +\tfrac{1}{2} \|x\|^2 \ ,
\end{equation}
and the first-order optimality conditions are simply:
$$
c+ (1+\|x\|_2^{r}) x = 0 \ ,
$$
whereby $x=-\theta c$ for some $\theta \ge 0$, and it remains to simply determine the value of the nonnegative scalar $\theta$.  If $c = 0$, then $x = 0$ satisfies the optimality conditions.  For $c \ne 0$, notice from above that $\theta$ must satisfy:
$$
1 - \theta - \|c\|_2^r \cdot \theta^{r+1}  =0 \ ,
$$
which is a univariate polynomial in $\theta$ with a unique positive root.  For $r =1,2,3$, this root can be computed in closed form.  Otherwise, the root can be computed (up to machine precision) using any scalar root-finding method. \medskip

\begin{rem} We can incorporate in problem \eqref{eq:rnsubprob} a simple set constraint $x\in Q$ provided that we can easily compute the Euclidean projection on $Q$. In the case when $\hh$ is a convex function of $\|x\|_2^2$, the subproblem \eqref{subproblem} can be converted to a $1$-dimensional convex optimization problem, see Appendix \ref{sec:constraintsubprob} for details.
\end{rem}  \medskip\medskip

{\bf A more specific example.}  Let $f(x) := \frac{1}{4} \| Ax- b\|_4^4 + \frac{1}{2} \|Cx-d\|_2^2$.  Then $\nabla^2 f(x)=3A^T D^2(x) A + C^T C$, where $D(x)= \BDiag(Ax-b)$.  Let us show that $ f(x)$ is $L$-smooth relative to  $$h(x) := \tfrac{1}{4} \|x\|_2^{4} + \tfrac{1}{2} \|x\|_2^2$$on $Q=\bbR^n$ for $L = 3\|A\|^4 + 6\|A\|^3 \|b\|_2 + 3\|A\|^2 \|b\|_2^2 + \|C\|^2$.  To see this, notice first that:\medskip
\begin{equation*}
\begin{array}{lcl}
\| \nabla^2 f(x) \| & \le & 3\|A\|^2 (\|b\|_2+\|A\|\|x\|_2)^2 + \|C\|^2 \\ \\
& = &\left( 3\|A\|^2 \|b\|_2^2 + \|C\|^2 \right) + 6\|A\|^3\|b\|_2 \|x\|_2 + 3\|A\|^4 \|x\|_2^2 \ , \\ \\
\end{array}
\end{equation*}\medskip which is $2$-degree polynomial in $\|x\|_2$ with coefficients $a_0 = 3\|A\|^2 \|b\|_2^2 + \|C\|^2$, $a_1 = 6\|A\|^3\|b\|_2$, and $a_2 = 3\|A\|^4$. Therefore following Remark \ref{uno} it suffices to set $$L = \sum_{i=0}^2 a_i = 3\|A\|^4 + 6\|A\|^3 \|b\|_2 + 3\|A\|^2 \|b\|_2^2 + \|C\|^2\ . $$

{\bf An example with Non-Lipschitz $\mu$-strong convexity.}  Let $f(x) := \frac{1}{4} \| Ex\|_2^4 + \frac{1}{4} \| Ax- b\|_4^4 + \frac{1}{2} \|Cx-d\|_2^2$, and let $\sigma_E$ and $\sigma_C$ denote the smallest singular values of $E$ and $C$, respectively, and let us suppose that $\sigma_E > 0 $ and $\sigma_C > 0$.  Then $\nabla^2 f(x)= \|Ex\|_2^2 E^T E + 2 E^T Exx^T E^T E + 3A^T D^2(x) A + C^T C$, where $D(x)= \BDiag(Ax-b)$. Let us show that $ f(x)$ is $L$-smooth and $\mu$-strongly convex relative to $$h(x) := \tfrac{1}{4} \|x\|_2^{4} + \tfrac{1}{2} \|x\|_2^2$$on $Q=\bbR^n$ for $L = 3 \|E\|^4 + 3\|A\|^4 + 6\|A\|^3 \|b\|_2 + 3\|A\|^2 \|b\|_2^2 + \|C\|^2$ and $\mu = \min \{ \frac{\sigma_E^4}{3}, \sigma_C^2 \}$.
Similar to what we have above,\medskip
\begin{equation*}
\begin{array}{lcl}
\| \nabla^2 f(x) \| & \le & \|E\|^4\|x\|_2^2 + 2\|E\|^4\|x\|_2^2+ 3\|A\|^2 (\|b\|_2+\|A\|\|x\|_2)^2 + \|C\|^2 \\ \\
& = &\left( 3\|A\|^2 \|b\|_2^2 + \|C\|^2 \right) + 6\|A\|^3\|b\|_2 \|x\|_2 + \left( 3 \|E\|^4 + 3\|A\|^4 \right) \|x\|_2^2 \ , \\ \\
\end{array}
\end{equation*}\medskip which is $2$-degree polynomial in $\|x\|_2$ with coefficients $a_0 = 3\|A\|^2 \|b\|_2^2 + \|C\|^2$, $a_1 = 6\|A\|^3\|b\|_2$, and $a_2 = 3 \|E\|^4 + 3\|A\|^4$. Therefore following Remark \ref{uno} it suffices to set $$L = \sum_{i=0}^2 a_i = 3 \|E\|^4 + 3\|A\|^4 + 6\|A\|^3 \|b\|_2 + 3\|A\|^2 \|b\|_2^2 + \|C\|^2\ . $$
On the other hand, \medskip
\begin{equation*}
\nabla^2 f(x) \succeq \|Ex\|_2^2 E^T E + C^T C \succeq \sigma_E^4 \|x\|_2^2 I + \sigma_C^2 I \succeq \mu \left( 1+3\|x\|_2^2 \right) I \succeq \mu \left(( 1+\|x\|_2^2)  I + 2 xx^T \right)= \mu \nabla^2 h(x)
\end{equation*}
(where the last matrix inequality follows since $\|x\|_2^2I \succeq xx^T$), and thus $f(x)$ is $\mu$-strongly convex relative to $h(x)$.\medskip\medskip

\begin{rem}\label{rem:center} In place of the simple reference function $\hh$ in \eqref{sunday} one can instead consider a ``re-centered'' version of the form:
\begin{equation*}\label{monday} h(x) = h_{x^c}(x) := \tfrac{1}{r+2} \|x-x^c\|_2^{r+2} + \tfrac{1}{2}\|x-x^c\|_2^2 \ , \end{equation*}
where the ``center'' value $x^c$ is suitably chosen to align $\ff$ with $\hh$ and possibly attain better values of $L$ and $\mu$.  Note that introducing the given center value $x^c$ does not increase the difficulty of solving the subproblem \eqref{subproblem}.  We illustrate this idea with a simple univariate example.  Suppose that our objective function is $f(x) = x^4 - 4x^3 +7x^2 - 5x + 3 $.  From the results in Section \ref{sec:rn} we know we can use the reference function $h_1(x) := \tfrac{1}{4}x^4 + \tfrac{1}{2}x^2$.  We can also translate $x$ by the center point $x^c:= 1$ and use the reference function $h_2(x):= \tfrac{1}{4}(x-1)^4 + \tfrac{1}{2}(x-1)^2$.  Straightforward calculation yields values of  $L=L_1=9+\sqrt{73} \approx 17.5440$ for $h_1(\cdot)$ and $L=L_2=4$ for $h_2(\cdot)$, whereby $h_2(\cdot)$ yields a better value of $L$ than $h_1(\cdot)$ for this example. \end{rem}

\subsection{$D$-Optimal Design Problem}\label{d-setup}

Given a matrix $H \in \mathbb{R}^{m \times n}$ of rank $m$ where $n \ge m + 1$, the $D$-optimal design problem is:
\begin{equation}\label{ddd}
\begin{array}{rll}
D: \ \ f^* = \min_x & f(x) := -\ln\det\left( H X H^T\right) \\ \\
\mathrm{s.t.} & \langle e, x \rangle = 1 \\
& x \ge 0 \ ,
\end{array}
\end{equation} where recall $X:= \BDiag(x)$.  In statistics, the $D$-optimal design problem corresponds to maximizing the determinant of the Fisher information matrix $\EE(hh^T)$, see  \cite{kiefer1960equivalence}, \cite{atwood1969optimal}.  And in computational geometry, $D$-optimal design arises as a Lagrangian dual problem of the minimum volume covering ellipsoid (MVCE) problem, which dates back at least 60 years to \cite{john}, see Todd \cite{toddminimum} for a modern treatment.  Indeed, \eqref{ddd} is useful in a variety of different application areas, for example, computational statistics \cite{croux2002location} and data mining \cite{knorr2001robust}.  In terms of algorithms for solving \eqref{ddd}, Khachiyan and Todd \cite{khachiyan1993complexity} proposed a theory-oriented scheme based on interior-point methods, see also Zhang \cite{zhang1998interior} as well as \cite{sun2004computation} for more practical treatments using interior-point methods.  Khachiyan \cite{k} later proposed and analyzed a first-order method (equivalent to the Frank-Wolfe method) to solve \eqref{ddd}, which led to other works along this line including Yildirim \cite{alperman} and Ahipasaoglu, Sun, and Todd \cite{suntodd}.  The complexity analysis in these papers is very specialized for the $D$-optimal design problem.  In contrast, we will show how the Primal Gradient Scheme (Algorithm \ref{proxgradscheme}) can be applied to the $D$-optimal design problem; furthermore, in Section \ref{return} we will apply the complexity analysis of Section \ref{algorithms} for the Primal Gradient Scheme to the set-up of $D$-optimal design, along with a comparison of our convergence guarantees with the guarantees from prior literature.

Notice that \eqref{ddd} is an instance of \eqref{poi1} with $Q=\Delta_n := \{x \in \mathbb{R}^n : \langle e, x \rangle = 1, \ x \ge 0\}$.  Although strictly speaking, $\ff$ in \eqref{ddd} is not defined everywhere on the relative boundary of $Q$ and hence does not have gradients or Hessians everywhere on the relative boundary of $Q$, this will not be of concern.  For $\ff$ in \eqref{ddd} let us choose the reference function $\hh$ to be the logarithmic barrier function, namely $$h(x) := -\sum_{j=1}^n \ln(x_j) \ , $$ which is defined on the positive orthant $\mathbb{R}^n_{++}$.  The following proposition states that $\ff$ is $1$-smooth relative to $\hh$.\medskip\medskip

\begin{prop}\label{d-smooth-prop} Suppose $f(x)=-\ln\det\left( H X H^T\right)$, where $X=\BDiag (x)$. Then $\ff$ is $1$-smooth relative to $h(x)=-\sum_{j=1}^n \ln(x_j)$ on $\mathbb{R}^n_{++}$. \qed \end{prop}

{\bf Proof:} The gradient of $\ff$ is $\nabla f(x)=\diag(-C)$ and the Hessian of $\ff$ is $ \nabla^2 f(x) = C \circ C $, where $C :=H^T (HXH^T)^{-1} H$.  Let $U=HX^{\frac{1}{2}}$; then $U^T (UU^T)^{-1} U \preceq I$ since the left side of this matrix inequality is a projection operator, whereby $X^{\frac{1}{2}}H^T(HXH^T)^{-1}HX^{\frac{1}{2}}\preceq I$.  Multiplying this matrix inequality on the left and right by $X^{-\frac{1}{2}}$ then shows that $C\preceq X^{-1}$ .
Therefore,
\begin{equation}\label{eq:dodhessian}
\nabla^2 f(x) = C \circ C \preceq C \circ X^{-1} \preceq X^{-1} \circ X^{-1} = X^{-2} = \nabla^2  h(x) \ ,
\end{equation}
where the first and the second matrix inequality above each follows from the fact that $C \preceq X^{-1} $ and the Hadamard product of two symmetric positive semidefinite matrices is also a symmetric positive semidefinite matrix.  The result then follows using property (a-iii) of Proposition \ref{thm:equivdef}. \qed

\noindent {\bf Solving the subproblem \eqref{subproblem}.}  Let us see how we can solve the subproblem \eqref{subproblem} for $Q$ and $\hh$ given above.  The subproblem \eqref{subproblem} can be written as
\begin{equation*}
\min_{x \in \Delta_n} \ \langle c, x \rangle -\sum_{j=1}^n \ln(x_j)  \ ,
\end{equation*}
and the first-order optimality conditions are simply:
$$
x > 0, \ \langle e, x \rangle = 1, \ \text{and} \ c - X^{-1}e = -\theta e \
$$
for some scalar multiplier $\theta$.  Given $\theta$, it then follows that $x_j = 1/(c_j + \theta)$ for $j=1, \ldots, n$, and it remains to simply determine the value of the scalar $\theta$.  Now notice that $\theta$ must satisfy:
\begin{equation}\label{smale}
d(\theta) := \sum_{j=1}^n \frac{1}{c_j + \theta} \ \ -1 \ = \ 0 \
\end{equation}
for some $\theta$ in the interval ${\cal F} := (-\min_j\{c_j\}, \infty)$.  Notice that $d(\cdot)$ is strictly decreasing on ${\cal F}$, and $d(\theta) \rightarrow +\infty$ as $\theta \searrow -\min_j\{c_j\}$ and $d(\theta) \rightarrow -1$ as $\theta \rightarrow \infty$, whereby \eqref{smale} has a unique solution in ${\cal F}$.  Furthermore, as suggested by results in Ye \cite{ye} or \cite{bf2}, one can use Newton's method (or any other suitable scalar solution-finding method) to efficiently compute the solution of \eqref{smale} (up to machine precision) on the interval ${\cal F}$ .\medskip\medskip

\subsection{Generalized Volumetric Function Optimization}

For a given integer parameter $p > 0$, let us also study optimization on the simplex of the following generalization of the volumetric barrier function:

\begin{equation}\label{vvv}
\begin{array}{cl}
\min_x & f_p(x) = \ln\det\left( H X^{-p} H^T\right) \\
\mathrm{s.t.} & \langle e, x \rangle = 1 \\
& x \ge 0 \ ,
\end{array}
\end{equation} where the integer $p$ is the parameter of the volumetric function $f_p(\cdot)$, and $H\in \mathbb{R}^{m\times n}$ is a rank-$m$ matrix where $n\ge m+1$.  Here the feasible region is $Q=\Delta_n$.  Note that $f_p(\cdot)$ is a convex function when $p \ge 0$ (and $f_p(\cdot)$ is a concave function when $p=-1$).

Similar to the $D$-optimal design problem, $f_p(\cdot)$ is not defined everywhere on the boundary of $\mathbb{R}^n_{+}$, but this will not be a concern. The reference function $\hh$ we choose is the logarithmic barrier function, namely $$h(x) := -\sum_{j=1}^n \ln(x_j) \ ,$$ which is defined on $\mathbb{R}^n_{++}$. The following proposition states that $f_p(\cdot)$ is $p(p+1)$-smooth relative to $\hh$.\medskip\medskip

\begin{prop}\label{d-smooth-prop-p} $f_p(\cdot)$ is $p(p+1)$-smooth relative to $h(x)=-\sum_{j=1}^n \ln(x_j)$ on $\mathbb{R}^n_{++}$.  \end{prop}

{\bf Proof:} By elementary calculus, the gradient of $f_p(\cdot)$ is $$\nabla f_p(x) = -p \cdot \diag \left(  X^{-1/2-p/2} C X^{-1/2-p/2}\right)\ ,$$ and the Hessian of $f_p(\cdot)$ is $$ \nabla^2 f_p(x) = p(p+1) \BMdiag \left( X^{-1-p/2} C X^{-1-p/2}  \right)-p^2 X^{-1-p/2} \left( C \circ C \right) X^{-1-p/2}\ , $$ where $C :=  H^T (HX^{-p}H^T)^{-1} H $, and $\BMdiag(M)$ denotes the diagonal matrix whose entries are the diagonal components of the matrix $M$. Let $U=HX^{-p/2}$; then $U^T (UU^T)^{-1} U \preceq I$ since the left side of this matrix inequality is a projection operator. Therefore each diagonal component of $U^T (UU^T)^{-1} U$ does not exceed $1$, whereby we have $\BMdiag \left( U^T(UU^T)^{-1}U \right) \preceq I$. Therefore, $$\begin{array}{lcl}
\nabla^2 f_p(x) & \preceq & p(p+1) \BMdiag \left( X^{-1-p/2} C X^{-1-p/2}  \right) \\ \\
& = & p(p+1) X^{-1} \BMdiag \left( U^T(UU^T)^{-1}U \right) X^{-1} \\ \\
& \preceq & p(p+1) X^{-2} \\ \\
& = & p(p+1) \nabla^2 h(x) \ ,
\end{array} $$
where the first inequality follows from the fact that the Hadamard product of two symmetric positive semidefinite matrices is also a symmetric positive semidefinite matrix and $C$ is a positive semidefinite matrix, and the first equation follows since $X$ is itself a diagonal matrix.  The result then follows by property (iii) of Proposition \ref{thm:equivdef}. \qed\medskip\medskip

\noindent {\bf Solving the subproblem \eqref{subproblem}.}  Using $h(x)=-\sum_{j=1}^n \ln(x_j)$, the subproblem \eqref{subproblem} here is identical to that for the $D$-optimal design problem, since the reference function $\hh$ and the feasible domain $Q$ are the same.  Therefore the methodology discussed in Section \ref{d-setup} applies here as well.  \medskip\medskip

{\bf Remark.} By setting $H=A^T$ and using Proposition \ref{prop:additivity}, it can also be shown that $\hat f(x):=\ln \det \left( A^T \BDiag \left( Ax-b \right)^{-p} A \right) $ is $p(p+1)$-smooth relative to $h(x):=-\sum_i \ln(A_i x -b_i)$. When $p=2$ this is the volumetric barrier function on the set $Q = \{x \in \mathbb{R}^n : Ax\ge b\}$, see \cite{vaidya1989new}, \cite{anstreicher2000volumetric}.

\subsection{Optimization over $Q \subset (0,u]^n$ with $\|\nabla^2f(x)\|$ growing as a polynomial in $\sum_{i=1}^n\frac{1}{x_i}$  }\label{sec:box}

Suppose that $\ff$ is a twice-differentiable convex function on $Q \subset (0,u]^n$ and that $\| \nabla^2 f(x) \| \le q_s \left(\sum_{i=1}^n \frac{1}{x_i} \right)$, where $q_s(\alpha)=\sum_{i=0}^s a_i \alpha^i$ is an $s$-degree polynomial in $\alpha$.  (Recall $\| \nabla^2 f(x) \|$ denotes the operator norm of $\nabla^2 f(x)$ with respect to the $\ell_2$-norm on $\mathbb{R}^n$.)   Let $$h(x) := \frac{u^3}{2(s+1)}\left( \sum_{i=1}^n \tfrac{1}{x_i}\right)^{s+1}\ . $$  Then the following proposition states that $\ff$ is $L$-smooth relative to $\hh$ for an easily computable value $L$.  This implies that no matter how fast $\nabla f(x)$ grows as $x$ approaches the open boundary of the region $(0,u]^n$, $\ff$ is smooth relative to the simple reference function $\hh$, even though $\nabla f(\cdot)$ need not exhibit uniform Lipschitz continuity on $Q$. \medskip\medskip

\begin{prop}\label{lem:smoothness-box}
Suppose $\ff$ is twice differentiable on $Q$ and satisfies $\|\nabla^2 f(x)\| \le q_s \left(\sum_{i=1}^n \frac{1}{x_i} \right)$ where $q_s(\alpha)$ is an $s$-degree polynomial in $\alpha$.  Let $L$ be such that $q_s(\alpha)\le L \alpha^s$ for all $\alpha\ge \tfrac{n}{u}$.  Then $\ff$ is $L$-smooth relative to $h(x) = \frac{u^3}{2(s+1)}(\sum_{i=1}^n \frac{1}{x_i})^{s+1}$.
\end{prop}
{\bf Proof:}  Let $X:=\BDiag(x)$, and it follows from elementary rules of differentiation that
\begin{equation}
 \nabla^2 h(x)  =  u^3 \left( \sum_{i=1}^n\tfrac{1}{x_i} \right)^s X^{-3} + \frac{u^3s}{2} \left( \sum_{i=1}^n \tfrac{1}{x_i} \right)^{s-1} X^{-2} ee^T X^{-2}\ . \end{equation}
Therefore
\begin{equation}
\nabla^2 h(x) \succeq  u^3 \left( \sum_{i=1}^n\tfrac{1}{x_i} \right)^s X^{-3}
 \succeq  \left( \sum_{i=1}^n \tfrac{1}{x_i} \right)^s I
 \succeq   \tfrac{1}{L} q_s\left(\sum_{i=1}^n\tfrac{1}{x_i}\right) I
 \succeq  \tfrac{1}{L} \nabla^2 f(x) \ ,
\end{equation}
where the second matrix inequality uses $u \ge x_i$ and the third matrix inequality is due to $\sum_{i=1}^n\tfrac{1}{x_i} \ge\sum_{i=1}^n \tfrac{1}{u} =\tfrac{n}{u}$. Therefore $\ff$ is $L$-smooth relative to $\hh$ by part (iii) of Proposition \ref{thm:equivdef}.\qed\medskip\medskip

\begin{rem} Suppose $q_s(\alpha)=\sum_{i=0}^s a_i \alpha^i$. In Proposition \ref{lem:smoothness-box}, one simple way to set $L$ is to use $L= \sum_{i=0}^s | a_i| \left(\tfrac{u}{n} \right)^{i-s}$.  This implies for $\alpha\ge \tfrac{n}{u}$ that
\begin{equation}
q_s(\alpha) \le \sum_{i=0}^s |a_i|\alpha^i \le \left( \sum_{i=0}^s |a_i| \left(\tfrac{u}{n} \right)^{i-s}\right) \alpha^s = L \alpha^s \ .
\end{equation}\end{rem}

\noindent {\bf Solving the subproblem \eqref{subproblem}.}  Let us see how we can solve the subproblem \eqref{subproblem} for this class of optimization problems.  After rescaling $c$ by $u^3/2$, the subproblem \eqref{subproblem} can be equivalently written as
\begin{equation}\label{pig}
\min_{x \in (0,u]^n} \ \langle c, x \rangle +  \tfrac{1}{s+1} \left( \sum_{i=1}^n \tfrac{1}{x_i} \right)^{s+1} \ .
\end{equation}
Let $\theta=\left( \sum_{i=1}^n \frac{1}{x_i} \right)^s$, then the optimality conditions for \eqref{pig} can be written as:
\begin{equation}\label{eq:first-order-box}
x_i \ = \ \left\{ \begin{array}{lcl}
 u & \mathrm{if} & c_i \le \frac{\theta}{u^2} \\ \\
\sqrt{\frac{\theta}{c_i}} & \mathrm{for} & c_i > \frac{\theta}{u^2} \ ,
\end{array}\right. \end{equation} for $i=1, \ldots, n$.  For a given $\theta >0$, define $x_i(\theta)$ using the above rule \eqref{eq:first-order-box}, and it remains to simply determine the value of the positive scalar $\theta$ in the interval ${\cal F} :=  [\left(\frac{n}{u}\right)^s, \infty)$ that satisfies \begin{equation}\label{ios}d(\theta) :=  \theta - \left( \sum_{i=1}^n \frac{1}{x_i(\theta)} \right)^s = 0 \ . \end{equation}

Notice that $d(\cdot)$ is strictly increasing on ${\cal F}$, and $d\left(\left(\frac{n}{u}\right)^s \right) \le 0$ (since $x_i(\theta) \le u$ for any $\theta$) and $d(\theta) \rightarrow \infty$ as $\theta \rightarrow  \infty$.  Therefore \eqref{ios} has a unique solution in ${\cal F}$, which can be solved with high accuracy using any suitable root-finding method, for example binary search combined with $1$-dimensional Newton's method.\medskip\medskip


\begin{rem}
In a sense, there are basically two ways that a twice-differentiable convex function can fail to have a uniformly Lipschitz gradient: (i) when the Hessian grows without limit as $\|x\| \rightarrow \infty$, and/or (ii) when the Hessian grows without limit as $x \rightarrow x^0 \in \partial Q$. Section \ref{sec:rn} has provided a mechanism for constructing a reference function $\hh$ for case (i) when the growth is polynomial, and Section \ref{sec:box} has provided such a mechanism for case (ii) when the growth is polynomial.  By utilizing the additivity and linear transformation properties of relative smoothness in Proposition \ref{prop:additivity}, it should be possible to construct suitable reference functions for many convex functions of interest.
\end{rem}






\section{Computational Analysis for the Primal Gradient Scheme and the Dual Averaging Scheme}\label{algorithms}

In this section we present computational guarantees for two algorithms: the Primal Gradient Scheme (Algorithm \ref{proxgradscheme}) as well as a Dual Averaging Scheme (Algorithm \ref{dualgradscheme}).

\subsection{Analysis of Primal Gradient Scheme (Algorithm \ref{proxgradscheme})}\label{sec:analysisprimal}

Our main result for the Primal Gradient Scheme is the following sublinear and linear convergence bounds.\medskip

%

\begin{thm}\label{thmprimalgd}
Consider the Primal Gradient Scheme (Algorithm \ref{proxgradscheme}). If $\ff$ is $L$-smooth and $\mu$-strongly convex relative to $\hh$ for some $L>0$ and $\mu \ge 0$, then for all $k \ge 1$ and $x\in Q$, sequence $\{f(x^k)\}$ is monotonically decreasing, and the following inequality holds:

\begin{equation}\label{joint}
\ \ \ \ \ \ \ f(x^k)  - f(x)   \ \ \leq \ \ \frac{\mu D_h(x,x^0)}{\left(1+\frac{\mu}{L-\mu} \right)^k-1} \ \ \le \ \  \frac{L-\mu}{k} D_h(x, x^0) \ ,
\end{equation}
where, in the case when $\mu = 0$, the middle expression is defined in the limit as $\mu \rightarrow 0^+$. \qed
\end{thm}\medskip\medskip

The first inequality in \eqref{joint} shows linear convergence when $\mu>0$; indeed, in this case it holds that \begin{equation}\label{fathersday}\frac{\mu D_h(x,x^0)}{\left(1+\frac{\mu}{L-\mu} \right)^k-1} \le  L\left( 1-\tfrac{\mu}{L}\right)^k D_h(x,x^0) \ . \end{equation}  (This inequality holds trivially for $k=1$, and induction on $k$ establishes the result for $k \ge 2$.)  Furthermore, when $k$ is large the $-1$  term in the denominator of the left-hand side can be ignored which yields the asymptotic bound $\mu\left( 1-\tfrac{\mu}{L}\right)^k D_h(x,x^0)$.  The second inequality in \eqref{joint} shows an $O(1/k)$ sublinear convergence rate. In particular, the convergence rate in \eqref{joint} is $\frac{L}{k} D_h(x, x^0)$ when $\mu = 0$.

Note that Algorithm \ref{proxgradscheme} herein and the NoLips algorithm in \cite{bbt} as well as algorithm PGA-$\cal B$ in \cite{zhouaunified} are structurally identical (they are all instantiations of the primal gradient methodology).  However, the step-size rule in \cite{bbt} as well as the complexity analysis in \cite{bbt} depends on a symmetry measure of $D_h(\cdot,\cdot)$, namely $\alpha := \min_{x,y \ne x} D_h(x,y)/D_h(y,x)$, whereas there is no such dependence here.  The instantiation of Algorithm \ref{proxgradscheme} in \cite{bbt} uses a smaller ``step-size'' of $(1+\alpha)/2L$ as opposed to $1/L$ in the update computation in Algorithm \ref{proxgradscheme} (since it must always hold that $\alpha \le 1$), and \cite{bbt} proves a computational guarantee of $f(x^k)  - f(x)  \le \frac{2L}{(1+\alpha) k} D_h(x, x^0)$.  The bound in Theorem \ref{thmprimalgd} is better than this symmetry-based bound, but only by a multiplicative constant factor $(1+\alpha)/2$ when $\mu=0$; it is of course far better (linear convergence rather than sublinear convergence) when $\mu >0$.\medskip

The proof of the bound in Theorem \ref{thmprimalgd} relies on the following standard Three-Point Property:\medskip\medskip

\begin{lem} {\bf (Three-Point Property of Tseng \cite{tseng})} Let $\phi(x)$ be a convex function, and let $\Dh(\cdot, \cdot)$ be the Bregman distance for $\hh$. For a given vector $z$, let
\begin{equation*}
z^+ := \arg\min_{x\in Q} \left\{ \phi(x) + \Dh(x,z) \right\} \ .
\end{equation*}
Then
\begin{equation*}
\phi(x) + \Dh(x,z) \ge \phi(z^+) + \Dh(z^+, z) + \Dh(x,z^+)\ \   for\ all \ x\in Q \ . \qed
\end{equation*}
\end{lem}\medskip\medskip

{\bf Proof of Theorem \ref{thmprimalgd}}:  Define a parameter sequence $$C_k:= \frac{1}{\sum_{i=1}^k\left(\frac{L}{L-\mu}\right)^i} \stackrel{(\cdot)}{=} \frac{\mu}{L\left(\left(1+\frac{\mu}{L-\mu}\right)^k - 1\right)} \ ,$$ where the second equality ``$(\cdot)$'' follows from elementary geometric series' analysis, and holds only when $\mu >0$.  In particular, $C_k = \frac{1}{k}$ if $\mu=0$.  For any $x \in Q$ and $i \ge 1$ we have:
\begin{equation}\label{in:fastconvergence}
\begin{array}{lcl}
f(x^{i}) & \le & f(x^{i-1}) + \langle \nabla f(x^{i-1}), x^{i}-x^{i-1} \rangle + L \Dh(x^{i}, x^{i-1}) \\ \\
& \le & f(x^{i-1}) + \langle \nabla f(x^{i-1}), x-x^{i-1} \rangle + L \Dh(x, x^{i-1}) -L \Dh(x,x^{i}) \\ \\
& \le & f(x) + (L-\mu) \Dh(x,x^{i-1}) - L\Dh(x,x^i) \ ,
\end{array}
\end{equation}
where the first inequality follows from the definition of $L$-smoothness relative to $\hh$, the second
inequality is due to the Three-Point Property with $\phi(x)=\tfrac{1}{L}\left\langle \nabla f(x^{i-1}),x-x^{i-1}\right\rangle $
and $z=x^{i-1}$, $z^{+}=x^{i}$, and the last inequality uses the $\mu$-strong convexity of $\ff$ relative to $\hh$, which implies $\langle \nabla f(x^{i-1}), x-x^{i-1} \rangle \le f(x) - f(x^{i-1}) -  \mu \Dh(x,x^{i-1})$.  Substituting $x=x^{i-1}$ in \eqref{in:fastconvergence} shows in particular that $f(x^i) \le f(x^{i-1})$ which proves monotonicity of the sequence $\{f(x^i)\}$.

It then follows using induction and \eqref{in:fastconvergence} that
\begin{equation}
\sum_{i=1}^k \left( \frac{L}{L-\mu}\right)^i f(x^i) \le  \sum_{i=1}^k \left( \frac{L}{L-\mu}\right)^i f(x) + L\Dh(x, x^0) - \left( \frac{L}{L-\mu}\right)^k L\Dh(x,x^k) \ .
\end{equation}

Using the monotonicity of $f(x^i)$ and the nonnegativity of $D_h(x,x^k)$, this implies that
\begin{equation}\label{toes}
 \left( \sum_{i=1}^k \left( \frac{L}{L-\mu}\right)^i\right) \left(f(x^k)-f(x)\right) \le L\Dh(x, x^0) - \left( \frac{L}{L-\mu}\right)^k L\Dh(x,x^k) \le L\Dh(x, x^0) \ .
\end{equation}
By substituting in the equality $$ \sum_{i=1}^k \left( \frac{L}{L-\mu}\right)^i = \frac{1}{C_k }  $$ in \eqref{toes} and rearranging, we obtain
\begin{equation}
\begin{array}{lcl}
f(x^k)  - f(x)   & \le &  C_k  L \Dh(x, x^0) = \displaystyle\frac{\mu D_h(x,x^0)}{\left(1+\frac{\mu}{L-\mu} \right)^k-1} \ .
\end{array}
\end{equation}

The proof of the second inequality in \eqref{joint} follows by noting that $\left(1+\frac{\mu}{L-\mu} \right)^k \ge 1+ \frac{k\mu}{L-\mu}$. \qed

\subsection{Dual Averaging Scheme and Analysis}\label{elbow}

Another algorithm for solving our optimization problem \eqref{poi1} is the Dual Averaging Scheme \cite{nesterovda}, which we present here in Algorithm \ref{dualgradscheme}.  Somewhat akin to the Primal Gradient Scheme, the update step in the Dual Averaging Scheme also requires the solution of a subproblem exactly of the form \eqref{subproblem}. Notice that we need the coefficient $\mu$ of strong convexity in order to implement Algorithm \ref{dualgradscheme}, in contrast to the Primal Gradient Scheme (Algorithm \ref{proxgradscheme}).  One can always conservatively set $\mu \gets 0$ in Algorithm \ref{dualgradscheme} if no reasonable lower bound on best value of $\mu$ is known.

\begin{algorithm}
\caption{Dual Averaging Scheme with reference function $h(\cdot)$}\label{dualgradscheme}
$ \ $
\begin{algorithmic}
\STATE {\bf Initialize.}  Let $L$, $\mu$ and $\hh$
satisfying
Definitions \ref{smooth} and \ref{strong} be given.  \\\
\ \ \ \ \ \ \ \ \ \ \ \ \ \ \ Let $x^0$ be the ``$\hh$-center'' of $Q$, namely $x^0 \gets \arg\min_{x\in Q} \{h(x)\}$, satisfying $h(x^0) = 0$.\\
\medskip
At iteration $k$ :\\
\STATE  {\bf Perform Updates.}  Compute $f(x^k)$, $\nabla
f(x^k)$ , $a_{k+1} = \frac{1}{L-\mu} \left(\frac{L}{L-\mu}\right)^k$, and \\\medskip
\ \ \ \ \ \ \ \ \ \ \ \ \ \
$x^{k+1} \gets \arg\min_{x \in Q}
\left\{ h(x) + \sum_{i=0}^k
a_{i+1}\left( f(x^i) + \langle \nabla f(x^i), x-x^i
\rangle + \mu D_h(x,x^i) \right)\right\}$ .\\\medskip
\end{algorithmic}
\end{algorithm}

We have the following result regarding computational guarantees for the Dual Averaging Scheme.\medskip\medskip

\begin{thm}\label{thmdualgd}
Consider the Dual Averaging Scheme (Algorithm
\ref{dualgradscheme}). If $\ff$ is $L$-smooth and
$\mu$-strongly convex relative to $\hh$ with $L > \mu$,
then for all $k \ge 1$ and $x\in Q$, the following
inequality holds:
\begin{equation}\label{eq-DRate}
\min_{i=1,\ldots,k} \{f(x^i)\} - f(x) \le \frac{\mu h(x)}{\left( 1 + \frac{\mu}{L-\mu}\right)^k - 1} \le \frac{L-\mu}{k}h(x)\ ,
\end{equation}
where in the case $\mu = 0$, the middle expression is
defined as the limits as $\mu \to 0^+$.
\end{thm}

Similar to the result in Theorem \ref{thmprimalgd}, the first inequality in \eqref{eq-DRate} shows linear convergence when $\mu >0$, since
\begin{equation}\label{mothersday}\frac{\mu h(x)}{\left(1+\frac{\mu}{L-\mu} \right)^k-1} \le  L\left( 1-\tfrac{\mu}{L}\right)^k h(x) \ ; \end{equation}this follows using identical logic as in \eqref{fathersday}.

{\bf Proof of Theorem \ref{thmdualgd}:} Define $\psi_{k}(x) := h(x) +
\sum\limits_{i=0}^{k-1} a_{i+1}\left( f(x^i) + \langle
\nabla f(x^i), x-x^i \rangle + \mu D_h(x, x^i) \right)$ for $k \geq 0$ and
$\psi_k^* := \min\limits_{x \in Q} \psi_k(x)$, whereby $x^k=\argmin_{x\in Q} \psi_k(x)$ and $\psi_k(x^k)=\psi_k^*$. It follows from the definition of relative strongly convexity (Definition \ref{strong}) that for any $x\in Q$:
\begin{equation}\label{eq-Up}
\psi_k^*  \leq  h(x) + A_k f(x) \ ,
\end{equation}
where $$A_k := \sum\limits_{i=0}^{k-1} a_{i+1} \stackrel{(\cdot)}{=}  \frac{1}{\mu} \left[ \left(1 + \frac{\mu}{L-\mu}\right)^k - 1 \right]$$ for all
$k \geq 0$, and where the second equality ``$(\cdot)$'' above follows from elementary geometric series' analysis and holds only when $\mu >0$; note that $A_k =  \frac{k}{L}$ when $\mu = 0$.

The function $\psi_k(\cdot)$ is a sum
of a linear function and the reference function $h(\cdot)$
multiplied by the coefficient $1+\mu A_k$. Therefore $(1+\mu A_k)h(\cdot)$ and $\psi_k(\cdot)$ define the same Bregman distance, whereby for any $x\in Q$ it holds that:
\begin{equation}\label{eq-Growth}
(1+\mu A_k) D_h(x, x^k) = D_{\psi_k}(x, x^k) = \psi_k(x) - \psi_k(x^k) - \langle\nabla \psi_k(x^k), x-x^k \rangle \le \psi_k(x) - \psi_k^* \ ,
\end{equation}
where the last inequality utilizes $\psi_k(x^k)=\psi_k^*$ as well as the first order optimality condition of $x^k=\argmin_{x\in Q} \psi_k(x)$.  Therefore:
\begin{equation*}
\begin{array}{lcl}
\psi_{k+1}^* & = & \psi_{k+1}(x^{k+1})\\ \\
& = & \psi_k(x^{k+1}) + a_{k+1}\left( f(x^k)
+ \langle \nabla f(x^k), x^{k+1}-x^k \rangle + \mu
D_h(x^{k+1},x^k) \right)\\ \\
& \ge & \psi_k^* + a_{k+1}\left( f(x^k) +
\langle \nabla f(x^k), x^{k+1}-x^k \rangle + \left( \mu +
\frac{1}{a_{k+1}} (1+\mu A_k)\right) D_h(x^{k+1},x^k)
\right) \ ,
\end{array}
\end{equation*}
where the last inequality uses \eqref{eq-Growth} with $x = x^{k+1}$.  Taking into account that $ \mu + \frac{1}{a_{k+1}} (1+\mu
A_k) = \frac{1 + \mu A_{k+1}}{a_{k+1}} = \frac{1}{a_{k+1}}
\left(\frac{L}{L-\mu}\right)^{k+1} = L$, and using the
relative smoothness of $f(\cdot)$ (Definition \ref{smooth}), we obtain: $$\psi^*_{k+1} \geq \psi_k^* + a_{k+1} f(x^{k+1})\ .$$  It then follows by induction that:
\begin{equation}\label{eq:almostdone}
\sum_{i=0}^{k-1} a_{i+1}f(x^{i+1}) \le \psi_{k}^*\le h(x) + A_k f(x) \ ,
\end{equation}
where the second inequality is from \eqref{eq-Up}.  The proof is completed by rearranging \eqref{eq:almostdone} and taking the minimum over $i$.
\qed

\subsection{On Optimization Problems with a Composite Function}\label{sec:composite}
Sometimes we are interested in solving the {\em composite} optimization problem \cite{nesterov2013}:
\begin{equation}\label{poi2}
\begin{array}{lrlr} P:  & f^* := \ \  \mbox{minimum}_x & f(x) + P(x)\\ \\
& \mbox{ s.t. } & x \in Q \ ,
\end{array}
\end{equation}
under the same assumptions on $\ff$ and $Q$ as in \eqref{poi1}, but now the objective function includes another function $P(\cdot)$ that is assumed to be convex but not necessarily differentiable, and for which the following subproblem is efficiently solvable:
\begin{equation}\label{subproblem2}
x_{\text{new}} \gets \arg\min_{x \in Q}\{\langle c, x \rangle + P(x) + h(x)\} \ ,
\end{equation} for any given $c$.  Under this assumption it is straightforward to show that Algorithm \ref{proxgradscheme} naturally extends to cover the case of the composite optimization problem \eqref{poi2} (see \cite{bbt} and \cite{zhouaunified}) and that the computational guarantee in Theorem \ref{thmprimalgd} extends to composite optimization as well.  (Indeed, when $\mu=0$ this extension is implied in principle from \cite{zhouaunified}.)  It turns out that one can actually view composite optimization as working with the objective function $\bar f(\cdot)$ that is $1$-smooth relative to the  reference function $\bar h(\cdot) :=L\hh + P(\cdot)$.  However, the definition of the reference function $\hh$ has been premised on $\hh$ being differentiable on $Q$, which might not hold for $\bar h(\cdot)$ as just defined.  This can all be taken care of by a suitable modification of the theory, see Appendix \ref{app:composite} for details.

\subsection{Questions: Accelerated Methods, Conjugate Duality, Choosing the Reference Function}

We have shown here in Section \ref{algorithms} that the computational guarantees of two standard first-order methods for smooth optimization -- the Primal Gradient Scheme and the Dual Averaging Scheme -- extend in precise ways to the case when $\ff$ is $L$-smooth relative to the reference function $\hh$.  The proof techniques used here suggest that very many other first-order algorithms for smooth optimization should extend similarly with analogous computational guarantees.  However, we have not been able to extend any accelerated methods, i.e., methods that attain an $O(1/k^2)$ convergence guarantee such as \cite{nest05smoothing}, \cite{nesterovBook}, \cite{tseng}, to the relatively smooth case.  One avenue for further research is to answer the question whether one can develop computational guarantees for an accelerated method in the case when $\ff$ is $L$-smooth relative to the reference function $\hh$?

Another question that arises concerns conjugate (duality) theory for the setting of relatively smooth convex functions.  One simple result in conjugate duality theory is that when $\ff$ is $L$-smooth (relative to $h(\cdot):=\tfrac{1}{2}\|\cdot\|^2$) the conjugate function $f^*(\cdot)$ is $1/L$-strongly convex (relative to $h^*(\cdot):=\tfrac{1}{2}\|\cdot\|_*^2$), see \cite{zalinescu}.  Is there a way to develop a more general conjugate duality theory that yields an analogous result when $\ff$ is $L$-smooth relative to a general convex function $\hh$?

A third question is how can we choose the reference function $\hh$ in order to lower the value of the bounds in Theorems \ref{thmprimalgd} and \ref{thmdualgd}?  Several ways to think about this question are discussed in Appendix \ref{app:chooseh}.

\section{$D$-Optimal Design Revisited:  Computational Guarantees using the Primal Gradient or Dual Averaging Scheme}\label{return}

Let us now apply the computational guarantees for the Primal Gradient Scheme (Theorem \ref{thmprimalgd}) and the Dual Averaging Scheme (Theorem \ref{thmdualgd}) to the $D$-optimal design optimization problem \eqref{ddd} discussed in Section \ref{d-setup}.  Recall from the exposition in Section \ref{d-setup} that $Q=\Delta_n$ and $f(x) = -\ln\det(HXH^T)$ is $1$-smooth relative to the logarithmic barrier function \begin{equation}\label{right} h(x) = -\sum_{j=1}^n \ln(x_j) \ , \end{equation}
and that the subproblem \eqref{subproblem} is efficiently solvable.  The following theorem presents a computational guarantee for using the Primal Gradient Scheme to approximately solve the $D$-optimal design optimization problem \eqref{ddd}.\medskip\medskip

\begin{thm}\label{thmDoptimalprimal} Consider using the Primal Gradient Scheme (Algorithm \ref{proxgradscheme}) with the reference function \eqref{right} to solve the $D$-optimal design problem \eqref{ddd} using the initial point $x^0 = \frac{1}{n} e$, and suppose that $\varepsilon \le f(x^0) - f^*$.  If $$k \ \ \ge \ \ \frac{2n\ln\left( \frac{2(f(x^0) - f^*)}{\varepsilon}\right)}{\varepsilon} \ , $$ then $ f(x^k) -f^* \le \varepsilon$.
\end{thm}

{\bf Proof:} Let $\delta = \frac{\varepsilon}{2(f(x^0) - f^*)}$. Then $\delta \le \frac{1}{2}$ since $\varepsilon \le f(x^0) - f^*$. Let $\hx:= (1-\delta) x^* + \delta x^0$. It follows from the convexity of $\ff$ that
$$
f(\hx) \le (1-\delta ) f^* + \delta f(x^0) \ ,
$$
whereby
\begin{equation}\label{eq:shrink}
 f(\hx)- f^* \le  \delta (f(x^0)- f^*) \ .
\end{equation}

Meanwhile,
\begin{equation}\label{eq:boundD}
\Dh(\hx, x^0) = h(\hx) - h(x^0) - \langle \nabla h(x^0), \hx - x^0 \rangle =h(\hx) - h(x^0) \le -n \ln \left(\tfrac{\delta}{n}\right) + n\ln\left( \tfrac{1}{n}\right) = n \ln(1/\delta) \ ,
\end{equation}
where the second equality uses $\nabla h(x^0)=-n \cdot e$ which then implies $\langle \nabla h(x^0), \hx - x^0 \rangle = 0$, and the inequality follows since $\hx\ge (\delta/n)e$.  Therefore, for $k$ satisfying the inequality in the statement of the theorem, we have:
\begin{equation}
\begin{array}{lcl}
f(x^k) - f^* & = & f(x^k) -f(\hx) + f(\hx)- f^* \\ \\
& \le & \displaystyle\frac{\Dh(\hx,x^0)}{k} + \delta (f(x^0)- f^*) \\ \\
& \le &  \displaystyle\frac{n \ln(1/\delta)}{k} + \frac{\varepsilon}{2} \\ \\
& \le &   \varepsilon\ ,
\end{array}
\end{equation}
where the first inequality follows from Theorem \ref{thmprimalgd} using $x = \hx$, as well as \eqref{eq:shrink}, the second inequality is from \eqref{eq:boundD} and the definition of $\delta$, and the third inequality follows since $k \ge [2n\ln (1/\delta)]/\varepsilon$. \qed\medskip\medskip

\begin{rem}\label{duo} For the Dual Averaging Scheme (Algorithm \ref{dualgradscheme}), one obtains the identical bound as in Theorem \ref{thmDoptimalprimal}.  This is proved by following virtually the same logic as above, except we use Theorem \ref{thmdualgd} which bounds the smallest optimality gap using $h(x) - h(x^0)$ instead of $\Dh(x, x^0)$.  However, it follows from \eqref{eq:boundD} that these two quantities are the same in this case.  Also, in the case of the Dual Averaging Scheme, the relevant final quantity of interest is $\min_{i=1, \ldots, k} f(x^i) - f^*$ instead of $f(x^k) - f^*$.
\end{rem}\medskip\medskip

It is instructive to compare the computational guarantees in Theorem \ref{thmDoptimalprimal}/Remark \ref{duo} to those of the Frank-Wolfe method applied to $D$-optimal design (first analyzed by Khachiyan \cite{k} and re-evaluated in \cite{suntodd} based in part on work by Yildirim \cite{alperman}).  Table \ref{comp-table} shows such a comparison, where absolute constants have been suppressed in order to highlight the dependencies on particular quantities of interest.  The second column of Table \ref{comp-table} compares the iteration bound of the methods using the starting point $x^0 = (1/n)e$, where we emphasize that $\varepsilon$ is the target optimality gap for the $D$-optimal design problem.  While it follows from observations in \cite{k} that $f(x^0) - f^* \le m \ln(n/m)$ for $x^0 = (1/n)e$, we do not show this in Table \ref{comp-table}, as we wish to highlight where the dependence on the initial iterate arises.  Examining the first column of Table \ref{comp-table}, note that the number of iterations of the Primal Gradient Scheme (or Dual Averaging Scheme) can be less than that of the Frank-Wolfe method, especially when $\varepsilon$ is not too small and when $n \ll m^2$.  However, as the second column of Table \ref{comp-table} shows, the Frank-Wolfe method requires only $mn$ operations per iteration in the worst -- i.e., dense matrix -- case, as it does a rank-$1$ update of a matrix inverse in the computation of $\nabla f(x^k)$), whereas the Primal Gradient Scheme (or Dual Averaging Scheme) requires $m^2n$ operations per iteration in the dense case (it must re-compute a matrix inverse in order to work with $\nabla f(x^k)$).  Therefore the total bound on operations of the Frank-Wolfe method (shown in the last column of Table \ref{comp-table}) is superior.

The bound for the Frank-Wolfe method applied to the $D$-optimal design problem is based on analysis that is uniquely designed for evaluating the $D$-optimal design problem, and is not part of the general theory for the Frank-Wolfe method (that we are aware of).  Even though the Primal Gradient Scheme and the Dual Averaging Scheme have inferior computational guarantees to the Frank-Wolfe method applied to the $D$-optimal design problem, they are the first (that we are aware of) first-order methods for which one has a general theory (Theorems \ref{thmprimalgd} and \ref{thmdualgd}) that can be meaningfully applied to yield computational guarantees for the $D$-optimal design problem.  We hope that this analysis will spur further interest in developing improved algorithms for $D$-optimal design and its dual problem -- the minimum volume enclosing ellipsoid problem.

\begin{table}[h]
\centering
\small{
\begin{tabular}{c||ccc}
 &   & Operations &  \\
 & Iteration  & Per Iteration & Total Operations \\
Method & Bound & (dense case) &  Bound\\ \hline
Frank-Wolfe Method & $m \ln(f(x^0)-f^*) + \displaystyle\frac{m^2}{\varepsilon}$ & $mn$ & $m^2n \ln(f(x^0)-f^*) + \displaystyle\frac{m^3n}{\varepsilon}$ \\ \\
Primal Gradient Scheme & $ \displaystyle\frac{n\ln(f(x^0)-f^*) }{\varepsilon} + \displaystyle\frac{n \ln\left(\frac{1}{\varepsilon}\right)}{\varepsilon}$ & $m^2n$ & $\displaystyle\frac{m^2n^2\ln(f(x^0)-f^*) }{\varepsilon} + \displaystyle\frac{m^2n^2 \ln\left(\frac{1}{\varepsilon}\right)}{\varepsilon}$ \\
or Dual Averaging Scheme \\

\end{tabular}}
\caption{Comparison of the order of computational guarantees for the Frank-Wolfe Method \cite{k}, \cite{suntodd} and the Primal Gradient and Dual Averaging Schemes (Theorem \ref{thmDoptimalprimal} and Remark \ref{duo}) for $D$-optimal design.  All constants have been suppressed in order to highlight the dependencies on particular quantities of interest.  It also follows from \cite{k} that $f(x^0) - f^* \le m \ln(n/m)$ for $x^0 = (1/n)e$, which can be inserted in the above bounds as well.}
\label{comp-table}
\end{table}

\section*{Acknowledgement}
The authors are grateful to the three referees for their comprehensive efforts and their suggestions on ways to improve the readability of the paper.

\appendix
\section{Appendix}

\subsection{Solving the subproblem \eqref{subproblem} when $h(x)$ is a convex function of $\|x\|_2^2$ and $Q$ has simple constraints} \label{sec:constraintsubprob}

We consider the following subproblem:
\begin{equation}\label{eq:subproblemconstraint}
\min_{x \in Q} \ \  \langle c, x \rangle + h(x) \ ,
\end{equation}
where $ h(x) = g(\|x\|_2^2)$ and $g(\cdot)$ is a (univariate) closed convex function of $\|x\|_2^2$.  Let  $y:=\|x\|_2^2$ and define ${\cal D}:=\{\|x\|_2^2 : x \in Q\} \subset \mathbb{R}$, which is the domain of $g(\cdot)$.  Let $g^*(\cdot)$ denote the conjugate function of $g(\cdot)$, namely
$$ g^*(t):= \sup_{y\in {\cal D}} \{ ty - g(y)\} \ , $$
whose domain we denote by ${\cal D}^*$.  Since $g(\cdot)$ is a convex function, we know from conjugacy theory \cite{avriel} that $g(y) = \sup_{t \in D^* }\{ty- g^*(t)\}$. Therefore \eqref{eq:subproblemconstraint} becomes
\begin{equation*}
\begin{array}{rcl}
\min_{x \in Q} \left\{ \langle c, x \rangle + g(\|x\|_2^2) \right\} & = &
 \min_{x \in Q} \left\{ \sup_{t\in D^*} \left\{ \langle c, x \rangle  + t\|x\|_2^2 - g^*(t)\right\} \right\}\\ \\
& = & \sup_{t\in D^*} \left\{ - g^*(t) + \min_{x \in Q}\left\{ \langle c, x \rangle  + t\|x\|_2^2 \right\} \right\} \ ,
\end{array}
\end{equation*} where the second equality above holds whenever the min and the sup operators can be exchanged (which is akin to strong duality).  Notice that $\min_{x \in Q} \{\langle c, x \rangle  + t\|x\|_2^2\}$ is a Euclidean projection problem.  Therefore the subproblem \eqref{subproblem} becomes a $1$-dimensional concave maximization problem if the Euclidean projection problem can be easily solved and one can conveniently form and work with the univariate convex conjugate function $g^*(\cdot)$.

\subsection{Extension to Composite Optimization}\label{app:composite}

Here we discuss some details of the extension of the ideas and results of this paper to composite optimization as described in Section \ref{sec:composite}, using the definitions $\bbf(\cdot) := \ff + P(\cdot)$, and $\bh(\cdot) = Lh(\cdot) + P(\cdot)$ as defined in Section \ref{sec:composite}.  Note that $\bar f(\cdot)$ and $\bar h(\cdot)$ are not necessarily differentiable on $Q$ since they include the function $P(\cdot)$. However, we can use the equivalent condition from (a-ii) of Proposition \ref{thm:equivdef} to define relative smoothness.  Let us now show how convergence results for the Primal Gradient Scheme still hold in this more general setting using an extension of the proof of Theorem \ref{thmprimalgd}.

Let $g_P(x) \in \partial P(x)$ be a specific subgradient of $P(\cdot)$ at $x$, and we will use the same subgradient of $P(\cdot)$ at $x$ when constructing a subgradient of $\bbf(\cdot)$ and/or $\bh(\cdot)$, namely $g_{\bbf}(x) := \nabla f(x) + g_P(x)$ and $g_{\bh}(x) := L\nabla h(x) + g_P(x)$. Then Algorithm \ref{proxgradscheme} has the following update:
\begin{equation}\label{eq:compositefunc}
\begin{array}{lcl}
x^{i+1} & = & \arg\min_{x\in Q}\{ \bar f(x^i) + \langle g_{\bbf}(x^i) , x-x^i \rangle + D_{\bh}(x, x^i) \} \\ \\
& = & \arg\min_{x\in Q}\{ \bar f(x^i) + \langle \nabla f(x^i) + g_P(x^i) , x-x^i \rangle  + D_{Lh}(x, x^i) + P(x) - P(x^i) - \langle  g_P(x^i) , x-x^i \rangle\}  \\ \\
& = & \arg\min_{x\in Q}\{ f(x^i) + \langle \nabla f(x^i), x-x^i \rangle + LD_{h}(x, x^i) + P(x)\} \ ,
\end{array}
\end{equation}
where in the third equality above the term involving $g_P(x^i)$ arising in $\partial \bar f(x^i)$ cancels out the corresponding term involving $g_P(x^i)$ arising in $\partial \bar h(x^i)$ as part of the expansion of $D_{\bh}(x , x^i)$.  There is therefore no actual need to compute $g_P(x^i) \in \partial P(x^i)$ in the update.   Indeed, this update \eqref{eq:compositefunc} corresponds exactly to the update in the NoLips algorithm \cite{bbt} (up to the step-size) and the PGA-$\cal B$ algorithm in \cite{zhouaunified} (up to the step-size) for composite optimization.

The proof of the computational guarantee in Theorem \ref{thmprimalgd} can be generalized directly to the composite optimization setting as follows.  Let us denote $$s_i(x):=\bar f(x^i) + \langle g_{\bbf}(x^i) , x-x^i \rangle + D_{\bh}(x, x^i) =   f(x^i) + \langle \nabla f(x^i), x-x^i \rangle + LD_{h}(x, x^i) + P(x)\ . $$  Notice that $x^{i+1}=\arg\min_{x\in Q} s_{i}(x)$; therefore from the first-order optimality conditions there is a subgradient $g_{s_{i}}(x^{i+1}) \in \partial s_{i}(x^{i+1})$ for which $\langle g_{s_{i}}(x^{i+1}), x-x^{i+1} \rangle \ge 0$ for all $x\in Q$.  From the additivity property of subgradients, we can write $g_{s_{i}}(x^{i+1}) = \nabla f(x^{i}) + L\nabla h(x^{i+1}) - L \nabla h(x^{i}) + \bar g$ for some $\bar g \in \partial P(x^{i+1})$, and let us assign $g_P(x^{i+1}) := \bar g = g_{s_{i}}(x^{i+1}) - \nabla f(x^{i}) - L\nabla h(x^{i+1}) + L\nabla h(x^{i})$, which then is used to define the subgradient $g_{\bbf}(x^{i+1})$, $g_{\bar h}(x^{i+1})$, and the Bregman distance $D_{\bh}(x, x^{i+1})$ in the proof.  Recall that the Primal Gradient Scheme does not rely on the choice of subgradient of $P(x^{i+1})$, thus the choice of $g_P(x^{i+1})$ is only used in the proof and it is well-defined.

Utilizing the above method for specifying the subgradients of $P(\cdot)$ at each of the iterates $x^i$ of the Primal Gradient Scheme, we can prove the following more specialized form of the Three Point Property which we can use in the proof of Theorem \ref{thmprimalgd} for the setting composite optimization.\medskip

\begin{lem} For any $x\in Q$, we have for any $i\ge 0$,
\begin{equation}\label{eq:g3p}
f(x^{i}) + \langle g_{\bbf}(x^{i}), x^{i+1}-x^{i} \rangle + D_{\bh}(x^{i+1}, x^{i}) \le f(x^{i}) + \langle g_{\bbf}(x^{i}), x-x^{i} \rangle + D_{\bh}(x, x^{i}) -  D_{\bh}(x, x^{i+1}) \ .
\end{equation}
\end{lem}
{\bf Proof:} Notice that $s_i(x)-\bh(x)= f(x^i) + \langle \nabla f(x^{i}) - L\nabla h(x^{i}) , x-x^{i} \rangle - Lh(x^{i})$ and so is a linear function of $x$, whereby it holds that
\begin{equation*}
\begin{array}{rcl}
(s_i(x) - \bh(x)) - (s_i(x^{i+1}) - \bh(x^{i+1}))
& = & \langle \nabla (s_i - \bh)(x^{i+1}), x-x^{i+1} \rangle \\ \\
& = & \langle g_{s_i}(x^{i+1}), x-x^{i+1} \rangle - \langle g_{\bh}(x^{i+1}), x -x^{i+1}\rangle \\ \\
& \ge & - \langle g_{\bh}(x^{i+1}), x -x^{i+1}\rangle \ ,
\end{array}
\end{equation*}
where the inequality follows from the choice of $g_{s_i}(x^{i+1})$. Rearranging the above and recalling the definition of $s_i(x)$ then completes the proof.\qed

The proof of Theorem \ref{thmprimalgd} in the setting of composite optimization follows directly by replacing $\hh$, $\nabla \hh$, $\ff$ and $\nabla \ff$ by $\bh(\cdot)$, $g_{\bh}(\cdot)$, $\bbf(\cdot)$ and $g_{\bar f}(\cdot)$, respectively, and utilizing \eqref{eq:g3p} to deduce the second inequality in \eqref{in:fastconvergence}.

\subsection{Criteria for choosing the reference function $\hh$}\label{app:chooseh}
One natural question is how can we choose $\hh$ in order to lower the value of the bound in Theorem \ref{thmprimalgd}?  Let us consider the simple case when $\ff$ is twice differentiable and is not strongly convex, namely $\mu = 0$, and $\ff$ attains its optimum at some point $x^*$.  Then the convergence bound \eqref{joint} can be re-written as:
\begin{equation*}
\begin{array}{lcl}
f(x^k) - f(x^*) & \le & \tfrac{1}{k}D_{Lh}(x^*, x^0) \\ \\
& = & \tfrac{1}{k}D_{f}(x^*, x^0) + \tfrac{1}{k} \left(\int_0^1 \int_0^t (x^* - x^0)^T \left[ \nabla^2(Lh-f) (x^0 + s(x^* - x^0) ) \right] (x^* - x^0) \ ds \ dt \right) \ ,
\end{array}
\end{equation*}
where $\nabla^2(Lh-f)(y)$ is the Hessian of the ``gap function'' $L\hh - \ff$ at the point $y\in Q$.  Notice that the first term above is fixed independent of the choice of $\hh$ and $L$.  It follows from Proposition \ref{thm:equivdef} that if $\ff$ is $L$-smooth relative to $\hh$ then $ \nabla^2(Lh-f)(y) \succeq 0$ for any $y\in \intQ$, whereby the second term above is always nonnegative. Since we do not know $x^*$ in most cases, in order to make the bound smaller we want the Hessian $\nabla^2(Lh-f)(y)$ to be smaller for all $y\in \intQ$.

There is a trade-off between how small the Hessian $\nabla^2(Lh-f)(y)$ is and how hard it will be to solve the subproblem \eqref{subproblem}. If we choose $L\hh = \ff$, the Hessian of the gap function is $0$, but solving the subproblem \eqref{subproblem} is as hard as solving the original problem \eqref{poi1}.  On the other hand, in standard gradient descent we use $\hh = \tfrac{1}{2}\|\cdot\|_2^2$ in which case the subproblem \eqref{subproblem} can be easily solved, while the Hessian of the gap function can be huge -- thus implying a poorer convergence bound.  There are a number of ways to try to manage this trade-off. For example, in gradient descent with preconditioning we can use $\hh = \tfrac{1}{2}\|\cdot\|_B^2 := \sqrt{\langle \cdot , B \cdot \rangle}$, where $B$ is a computationally-friendly positive definite matrix -- typically a diagonal matrix. The criteria for designing $B$ usually involves (i) ensuring that solving equations with $B$ is easy (so that the subproblem \eqref{subproblem} can be easily solved), and (ii) $B$ is ``close to'' the Hessian of $\ff$ (so that the Hessian of the gap function is small).

\bibliographystyle{amsplain}
\bibliography{GF-papers-nips-better}
\end{document}